\documentclass[12pt,letterpaper]{amsart}

\newcounter{myenum}

\newcommand{\indentalign}{\hspace{0.3in}&\hspace{-0.3in}}

\newcommand{\R}{\operatorname{Re}}
\newcommand{\I}{\operatorname{Im}}
\newcommand{\ds}{\displaystyle}
\newcommand{\sech}{\operatorname{sech}}

\usepackage{amssymb}
\usepackage{amsthm}
\usepackage{amsxtra}

\newtheorem{theorem}{Theorem}

\newtheorem{proposition}[theorem]{Proposition}

\newtheorem{corollary}[theorem]{Corollary}

\theoremstyle{remark}
\newtheorem{remark}[theorem]{Remark}

\newcommand{\cR}{\mathbb{R}}

\newcommand{\ep}{\epsilon}

\newcommand{\lap}{\triangle}

\numberwithin{equation}{section}
\numberwithin{theorem}{section}

\title[Blow-up solutions to 3D NLS]
{On blow-up solutions to the 3D cubic nonlinear Schr\"odinger equation}
\author{Justin Holmer}
\address{University of California, Berkeley}
\author{Svetlana Roudenko}
\address{Arizona State University}

\thanks{to appear in {\it Appl. Math. Res. eXpress}\\
submitted Sept 28, 2006,\\
accepted Feb 12, 2007}

\begin{document}

\begin{abstract}
For the 3d cubic nonlinear Schr\"odinger (NLS) equation, which has
critical (scaling) norms $L^3$ and $\dot H^{1/2}$, we first prove a result establishing
sufficient conditions for global existence and sufficient conditions
for finite-time blow-up.  For the rest of the paper, we focus on the
study of finite-time radial blow-up solutions, and prove a result on the
concentration of the $L^3$ norm at the origin. Two disparate
possibilities emerge, one which coincides with solutions typically
observed in numerical experiments that consist of a specific bump
profile with maximum at the origin and focus toward the origin at
rate $\sim(T-t)^{1/2}$, where $T>0$ is the blow-up time.  For the
other possibility, we propose the existence of ``contracting sphere
blow-up solutions'', i.e. those that concentrate on a sphere of
radius $\sim (T-t)^{1/3}$, but focus towards this sphere at a faster
rate $\sim (T-t)^{2/3}$. These conjectured solutions are analyzed
through heuristic arguments and shown (at this level of precision)
to be consistent with all conservation laws of the equation.
\end{abstract}

\maketitle

\section{Introduction}
\label{S:intro}

Consider the cubic focusing nonlinear Schr\"odinger (NLS) equation on $\mathbb{R}^3$:
\begin{equation}
 \label{E:NLS}
i\partial_t u +\Delta u + |u|^2u=0,
\end{equation}
where $u=u(x,t)$ is complex-valued and $(x,t)\in
\mathbb{R}^3\times \mathbb{R}$. The initial-value problem posed with
initial-data $u(x,0)=u_0(x)$ is locally well-posed in $H^1$. Such
solutions, during their lifespan $[0,T)$ (where $T=+\infty$ or
$T<+\infty$), satisfy mass conservation $M[u](t)=M[u_0]$, where
$$
M[u](t) = \int |u(x,t)|^2 \,dx,
$$
and energy conservation $E[u](t)=E[u_0]$, where
$$
E[u](t) = \frac12\int |\nabla u(x,t)|^2 - \frac14\int |u(x,t)|^4
\,dx.
$$
If $\|xu_0\|_{L^2}<\infty$, then $u$ satisfies the virial identity
$$
\partial_t^2 \int |x|^2 |u(x,t)|^2 \, dx = 24E[u] - 4\|\nabla
u(t)\|_{L_x^2}^2.
$$
The equation has scaling: $u_\lambda(x,t) = \lambda u(\lambda x,
\lambda^2 t)$ is a solution if $u(x,t)$ is a solution. The
scale-invariant Lebesgue norm for this equation is $L^3$, and the
scale-invariant Sobolev norm is $\dot H^{1/2}$.

Fundamental questions include:
\begin{enumerate}
\item[1.]
Under what conditions on the initial data $u_0$ is the solution $u$
globally defined ($T=+\infty$)? If it is globally defined, does it
scatter (approach the solution to a linear Schr\"odinger equation as
$t\to +\infty$) or resolve into a sum of decoupled solitons plus a
dispersive component? The latter type of inquiry leads one to the
``soliton resolution conjecture'' (see Tao \cite{T04, T06}).

\item[2.]
If the solution fails to be globally defined (we say ``blows-up in
finite time''), can one provide a description of the behavior of the
solution as $t\to T$, where $T$ is the ``blow-up time''?  This will
be the focus of the present paper.
\end{enumerate}
It follows from the $H^1$ local theory optimized by scaling (see
Cazenave \cite{Caz-book} or Tao \cite{Tao-book} for exposition),
that if blow-up in finite-time $T>0$ occurs, then there is a
lower-bound on the ``blow-up rate'':
\begin{equation}
\label{E:lowerbd}
\|\nabla u(t)\|_{L_x^2} \geq \frac{c}{(T-t)^{1/4}}
\end{equation}
for some absolute constant\footnote{By ``absolute constant'' in this
article we mean a constant that does not depend on any properties of
the solution $u$ under study (e.g.\ mass, energy, etc.), and usually
depends on constants appearing in Gagliardo-Nirenberg, Sobolev, and
Strichartz estimates.} $c$. Thus, to prove global existence, it
suffices to prove a global \textit{a priori} bound on $\|\nabla
u(t)\|_{L^2}$. From the Strichartz estimates, there is a constant
$c_{ST}>0$ such that if $\|u_0\|_{\dot H^{1/2}}< c_{ST}$, then the
solution $u$ is globally defined and scatters.  The
\textit{optimal}\footnote{This is not the ``sharp" constant in
Strichartz inequalities but rather the one which will govern
scattering.} constant $c_{ST}$, the ``scattering threshold'', has not
to our knowledge been identified, under any regularity or decay
assumptions on the solution.

Note that the quantities $\|u_0\|_{L^2}\|\nabla u_0\|_{L^2}$ and
$M[u_0]E[u_0]$ are also scale-invariant.  Another result of the
above type follows from the conservation laws and the sharp
Gagliardo-Nirenberg inequality of M. Weinstein \cite{W83}.  Let
$Q(x)$ be the minimal mass ground-state solution to the nonlinear
elliptic equation
$$-\tfrac12Q+\tfrac32 \Delta Q+Q^3=0$$
on $\mathbb{R}^3$, and set $u_Q(x,t) = e^{\frac12 it}Q(\sqrt\frac32
x)$. Then $u_Q$ is a soliton solution to \eqref{E:NLS}, and we have:

\begin{theorem}
\label{T:glob_vs_blowup} Suppose $u$ is the (possibly nonradial)
$H^1$ solution to \eqref{E:NLS} corresponding to initial data $u_0$
satisfying
\begin{equation}
\label{E:mass_energy_threshold}
M[u_0]E[u_0] < M[u_Q]E[u_Q] \, .
\end{equation}

{\rm 1.} If $\|u_0\|_{L^2}\|\nabla u_0\|_{L^2}< \|u_Q\|_{L^2}\|\nabla u_Q\|_{L^2}$, then
$\|u_0\|_{L^2}\|\nabla u(t)\|_{L^2}< \|u_Q\|_{L^2}\|\nabla u_Q\|_{L^2}$ for all time, and
thus, the solution is global.\footnote{Although, as far as we are
aware, the threshold for scattering \textit{might} be strictly
smaller than  $\|u_Q\|_{L^2}\|\nabla u_Q\|_{L^2}$. } %

{\rm 2.} %
If $\|u_0\|_{L^2}\|\nabla u_0\|_{L^2}> \|u_Q\|_{L^2}\|\nabla
u_Q\|_{L^2}$, then $\|u_0\|_{L^2}\|\nabla u(t)\|_{L^2}>
\|u_Q\|_{L^2}\|\nabla u_Q\|_{L^2}$ on the maximal time interval of
existence. If we further assume finite variance
$\|xu_0\|_{L^2}<\infty$ or radial symmetry of the solution, then the
solution blows-up in finite time.

If \eqref{E:mass_energy_threshold} fails to hold, then we have no
information to conclude global existence or finite-time blow-up.
\end{theorem}
We note that if the energy is negative, then via the sharp
Gagliardo-Nirenberg inequality we automatically have
$\|u_0\|_{L^2}\|\nabla u_0\|_{L^2}>\|u_Q\|_{L^2}\|\nabla u_Q\|_{L^2}$, and the second of
the two cases in the theorem applies.  In \S \ref{S:finite-time}, we
prove a generalized version of Theorem \ref{T:glob_vs_blowup}, and
explain that it provides a bridge between similar known results for
the $L^2$-critical (by Weinstein \cite{W83}) and $\dot H^1$-critical
(by Kenig-Merle \cite{KM}) NLS equations.

Now, suppose blow-up in finite time occurs (for a solution of any
energy), and let us restrict attention to the radially symmetric
case.  What can be said about the behavior of the solution as $t\to
T$, where $T$ is the blow-up time?  As mentioned above, there is a
lower bound \eqref{E:lowerbd} on the blow-up rate. Merle-Rapha\"el
\cite{MR06} have recently shown that the scale-invariant norm
$\dot{H}^{1/2}$ has a divergent lower-bound:
$$
\|u(\cdot,t)\|_{\dot{H}^{1/2}} \geq c\, |\log (T-t)|^{1/12}.
$$
This is in contrast to the $L^2$-critical problem, where the scale
invariant norm is, of course, constant.  Merle-Rapha\"el do not
prove any upper bound for general solutions on the rate of
divergence of this norm.\footnote{Indeed, the heuristic analysis of
contracting sphere blow-up solutions we provide in this article
suggests the existence of a solution for which
$\|u(\cdot,t)\|_{\dot{H}^{1/2}} \sim (T-t)^{-1/3}$, far larger than
the general lower bound of \cite{MR06}.  Merle-Rapha\"el do show,
however, that if equality is achieved in \eqref{E:lowerbd}, one can
obtain an upper bound $\|u(t)\|_{\dot H^{1/2}} \leq |\log
(T-t)|^{3/4}$.} The second result we present in this note describes
two possibilities for the rate of concentration of the $L^3$ norm
for radial finite-time blow-up solutions.  We find it more
convenient in this analysis to work with the critical Lebesgue norm
$L^3$, rather than the critical Sobolev norm $\dot H^{1/2}$, since
the former is more easily localized.   In the physics or numerics
literature (see Sulem-Sulem \cite{SS99}), this type of concentration
phenomena is termed ``weak concentration'' to distinguish it from
concentration in the $L^2$ norm, which is called ``strong
concentration''.

\begin{theorem}
\label{T:L3conc} Suppose the radial $H^1$ solution $u$ to
\eqref{E:NLS} blows-up at time $T<\infty$. Then either there is a
non-absolute\footnote{This means it depends on the solution but not
on time.} constant $c_1 \gg 1$ such that, as $t\to T$
\begin{equation}
\label{E:tight}
\int_{|x|\leq c_1^2\|\nabla u(t)\|_{L^2}^{-2}} |u(x,t)|^3 \, dx\geq c_1^{-1}
\end{equation}
or there exists a sequence of times $t_n\to T$ such that, for an
absolute constant $c_2$
\begin{equation}
 \label{E:wide}
\int_{|x|\leq c_2\|u_0\|_{L^2}^{3/2}\|\nabla u(t)\|_{L^2}^{-1/2}}
|u(x,t_n)|^3 \, dx \to \infty \,.
\end{equation}
\end{theorem}
These two cases are not mutually exclusive. From the lower bound
\eqref{E:lowerbd}, we have that the concentration window in
\eqref{E:tight} satisfies $\|\nabla u(t)\|_{L^2}^{-2} \leq
c\,(T-t)^{1/2}$, and the concentration window in \eqref{E:wide}
satisfies $\|\nabla u(t)\|_{L^2}^{-1/2} \leq c\,(T-t)^{1/8}$.  The
argument combines the radial Gagliardo-Nirenberg estimate (as we
learned from J.\ Colliander, private communication) and the argument
in the proof of Proposition 7 in Hmidi-Keraani \cite{HK05}.  It may
be that a more refined analysis using the successive extraction of
weak limits technique in %
\cite{HK05} could yield more
precise information,
although we have decided not to explore this for the moment.  This
result can be compared with the mass concentration result for $L^2$
critical equations, see Merle-Tsutsumi \cite{MT90} for the first
results in this direction (radial case) and the recent paper of
Hmidi-Keraani \cite{HK05} for references and a simplified proof in
the general (nonradial) case.

This analysis led us to consider: What type of blow-up solution
would display the behavior described in \eqref{E:tight}, and what
type would display the behavior described in \eqref{E:wide}? There
are currently no analytical results describing the specifics of the
profile of the solution as $t\to T$ for finite-time blow-up
solutions, although there have been several numerical studies.  We
mention, however, the construction by Rapha\"el \cite{R06} of a
solution to the two-dimensional quintic NLS (also
mass-supercritical) that blows-up on the unit circle. Rapha\"el's
result \cite{R06} draws upon a large body of breakthrough work by
Merle-Rapha\"el \cite{MR-Ann}-\cite{MR06} \cite{R06} (see also
Fibich-Merle-Rapha\"el \cite{FMR06}) on the blow-up problem for the
$L^2$ critical NLS. Rapha\"el's construction \cite{R06}, and the
numerical study of Fibich-Gavish-Wang \cite{FGW05}, have inspired
our inquiry into the ``contracting sphere'' solutions that we
describe through a heuristic analysis below. First, however, we call
attention to the numerical results (see Sulem-Sulem \cite{SS99} for
references) describing the existence of self-similar radial blow-up
solutions of the form
\begin{equation}
 \label{E:atorigin}
u(x,t) \approx \frac{1}{\lambda(t)}U\left( \frac{x}{\lambda(t)}
\right) e^{i\log(T-t)} \text{ with }\lambda(t)=\sqrt{2b(T-t)}
\end{equation}
for some parameter $b>0$ and some stationary profile $U=U(x)$
satisfying the nonlinear elliptic equation
$$
\Delta U - U + i b ( U + y\cdot\nabla U) + |U|^2 U =0.
$$
It is expected that (nontrivial) zero-energy solutions $U$ to this
equation exist, although fail to belong to $L^3$ (and thus also
$\dot H^{1/2}$) due to a logarithmic growth at
infinity.\footnote{The only rigorous results on the existence of
such $U$ are for mass-supercritical NLS equations that scale in
$\dot H^{s_c}$, for $s_c>0$ close to $0$ -- see Rottsch\"afer-Kaper
\cite{RK}.} Thus, for solutions in $H^1$, the interpretation of
$\approx$ in \eqref{E:atorigin} is that one should introduce a
time-dependent truncation of the profile $u$, where the size of the
truncation enlarges as $t\to T$. The resulting $u$ described by
\eqref{E:atorigin} will then display at least the logarithmic $\dot
H^{1/2}$ divergence that must necessarily occur by the work of
Merle-Rapha\"el \cite{MR06}.\footnote{We thank J.\ Colliander for
supplying in private communication this interpretation derived from
discussions with C.\ Sulem.} We note that this type of solution
would also display the $L^3$ concentration properties described in
\eqref{E:tight} in Theorem \ref{T:L3conc}, and in fact, would also
satisfy \eqref{E:wide}.

We are then led to consider: Could there be a solution for which
\eqref{E:wide} holds but \eqref{E:tight} does not hold?  Such a
solution would have to concentrate on a contracting sphere of radius
$r_0(t)$, where $r_0(t)\to 0$. That is, essentially no part of the
solution sits directly on top of the origin $|x|\leq \frac12
r_0(t)$, and the whole blow-up action is taking place inside the
spherical annulus $\frac12r_0(t)<|x|<\frac32r_0(t)$.  Conservation
of mass dictates that in fact the rate of contraction of the
solution towards the sphere must far exceed the rate of contraction
of the sphere itself.  Specifically we seek a solution of the form
(in terms of amplitudes)
$$
|u(r,t)| \approx \frac{1}{\lambda(t)} \, P\left(
\frac{r-r_0(t)}{\lambda(t)}\right),
$$
where $P$ is some one-dimensional profile and $r$ is the
three-dimensional radial coordinate. Then, we must have $\lambda(t)
\sim r_0(t)^2$.  By studying all the conservation laws  and allowing
for a little more generality, we provide a heuristic argument
suggesting that such solutions do exist, but only with the following
specific features.

\vspace{.3cm}

\noindent\textbf{Blow-up Scenario}.  \textit{With blow-up time $0 < T < \infty$,
we define for $t<T$ the
radial position and focusing factor
$$
r_0(t) = \frac{3^{1/6}M[u]^{1/3}}{2\pi^{1/3}}(T-t)^{1/3}, \qquad
\lambda(t) = \frac{18^{1/3}}{M[u]^{1/3}}(T-t)^{2/3}
$$
and the rescaled time
$$s(t) = \frac{3M[u]^{2/3}}{18^{2/3}}(T-t)^{-1/3}$$
(so that, in particular, $s(t)\to +\infty$ as $t\to T$). Label the
constant $\kappa = \frac43\left( \frac{3}{32\pi}\right)^{2/3}$, and
take $\theta$ as an arbitrary phase shift.  Then
\begin{equation}
 \label{E:sphere_blowup}
u(r,t) \approx e^{i\theta}e^{i\kappa^2s}\exp\left[ i\kappa
\left(\frac{r-r_0(t)}{2\,\lambda(t)}\right)\right] \,
\frac1{\lambda(t)} \, P\left( \frac{r-r_0(t)}{\lambda(t)}\right),
\end{equation}
where
$$
P(y) = \sqrt\frac32 \,\kappa \, \sech\left(\frac{\sqrt 3}{2} \,
\kappa y\right),
$$
is a blow-up scenario that is consistent with all
conservation laws.  }

\vspace{.3cm}

A crucial component of this analysis is the
observation that the inclusion of the spatial phase-shift gives a
profile of zero energy, which is essential since the rescaling of
the solution through the focusing factor $\lambda(t)$ would cause a
nonzero energy to diverge to $+\infty$.  We believe, however, that
it is possible for the actual solution $u$ to have nonzero energy,
since it will in fact be represented in the form $u(x,t) =
\text{(above profile)}+\tilde u(x,t)$, where $\tilde u(x,t)$ is an
error, and $\tilde u(x,t)$ could introduce nonzero energy by itself
or through interaction with the profile \eqref{E:sphere_blowup}.

The analysis demonstrates that the mass conservation and the virial
identity act in conjunction to drive the blow-up:  start with an
initial configuration of the form \eqref{E:sphere_blowup}, then the
virial identity will tend to push $r_0(t)$ inward, and the mass
conservation will force the solution to compensate by driving
$\lambda(t)$ smaller (focusing the solution further), which in turn
will feed back into the virial identity to push $r_0(t)$ yet
smaller. We emphasize that we do not actually prove here that such
blow-up solutions exist (we are currently working on such a rigorous
construction which will be presented elsewhere), we only provide
evidence through heuristic reasoning in this note. The numerical
observation of solutions of this type was reported in
Fibich-Gavish-Wang \cite{FGW05}, although most of the specifics of
the dynamic were not mentioned -- only the relationship
$\lambda(t)\sim r_0(t)^2$.  The authors did remark that they
intended to follow-up with a more detailed analysis of this
supercritical problem, and devoted most of the \cite{FGW05} paper to
the analysis of a similar phenomenon in the $L^2$ critical context.

For the hypothetical contracting sphere solutions, we have, as $t\to T$,
\begin{equation}
\label{E:rates}
\begin{aligned}
&\|u(t)\|_{L^3} \sim \frac{r_0(t)^{2/3}}{\lambda(t)^{2/3}}
\sim M[u_0]^{4/9}(T-t)^{-2/9}, \\
&\|u(t)\|_{\dot H^{1/2}} \sim \frac{r_0(t)}{\lambda(t)}
\sim M[u_0]^{2/3}(T-t)^{-1/3},\\
&\|\nabla u(t)\|_{L^2} \sim \frac{r_0(t)}{\lambda(t)^{3/2}} \sim
M[u_0]^{5/6}(T-t)^{-2/3}.
\end{aligned}
\end{equation}
Thus, the concentration window in \eqref{E:wide} is $\sim
(T-t)^{1/3}$ (which coincides with $r_0(t)$) and in \eqref{E:tight}
is $\sim (T-t)^{4/3}$.  This type of solution, if it exists, should
satisfy \eqref{E:wide} without satisfying \eqref{E:tight}.  Another
remark in regard to \eqref{E:rates} is that \eqref{E:sphere_blowup}
is only an asymptotic description for $t$ close to $T$, where this
closeness depends on $M[u_0]$ (and potentially other factors), and
thus, there is no contradiction with the above formula for
$\|u(t)\|_{\dot H^{1/2}}$ and the small data $\dot H^{1/2}$
scattering theory.  A similar comment applies to the other scale
invariant quantity $\|u(t)\|_{L^2}\|\nabla u(t)\|_{L^2}$ discussed
earlier.  (In addition, there is no clear constraint on the energy
of such solutions -- they may satisfy
\eqref{E:mass_energy_threshold} or may not.)

Since the contracting sphere blow-up solutions contract at the rate
$\sim (T-t)^{1/3}$, and the solutions \eqref{E:atorigin} that
blow-up on top of the origin focus at a faster rate $\sim
(T-t)^{1/2}$, we propose the possibility of a solution that blows-up
simultaneously in both manners -- we see from the rates that a
decoupling should occur between the two components of the solution.
One could also speculate that multi-contracting sphere solutions,
with or without a contracting blob at the origin itself, are
possible.  We further conjecture, on the basis of the article
\cite{FGW05}, that the contracting sphere blow-up solutions are
stable under small radial perturbations but unstable under small
nonradial (symmetry breaking) perturbations\footnote{Stability under
small $H^1$ radial perturbations was also proved for the singularity
formation on a constant ring in \cite{R06}, and it was asked there whether the
non-radial perturbations will affect the stability.}.

We remark that for the defocusing nonlinearity ($+|u|^2u$ changed to
$-|u|^2u$ in \eqref{E:NLS}), one always has global existence, since
the energy is then a positive definite quantity and thus
automatically provides an \textit{a priori} bound on the $\dot H^1$
norm.  Furthermore, it has been shown by Ginibre-Velo \cite{GV85},
using a Morawetz estimate, that there is scattering.  The argument
of Ginibre-Velo has been simplified using a new ``interaction
Morawetz'' identity by Colliander-Keel-Staffilani-Takaoka-Tao
\cite{CKSTT04}.  Thus, at least as far as $H^1$ data is concerned,
the dynamics of this problem are comparatively well-understood.

The format for this note is as follows. In \S \ref{S:finite-time},
we study a general version of the focusing NLS equation which is
energy subcritical and prove the generalized version of Theorem
\ref{T:glob_vs_blowup}; there we also discuss the blow up criterion
for $u_0 \in H^1$ which includes positive energies and not
necessarily finite variance. In \S \ref{S:conc}, we prove Theorem
\ref{T:L3conc} on $L^3$ norm concentration.  In \S
\ref{S:heuristic}, we present the heuristic analysis of the
conjectured contracting sphere solutions. This is followed in \S
\ref{S:refined} by an analysis with somewhat more precision,
indicating, in particular, a cancelation between second-order
approximations to two different ``error'' terms.  We conclude in \S
\ref{S:general} by noting that the ideas presented here for the
specific problem \eqref{E:NLS} can be adapted to more general radial
nonlinear Schr\"odinger equations. Interestingly, there is the
possibility that for septic nonlinearity $|u|^6 u$, one could have
blow-up on an \textit{expanding} sphere, but we have not conducted a
thorough analysis of this hypothetical situation.\\

\noindent\textbf{Acknowledgments}. We would like to thank the
organizers of the Fall 2005 MSRI program ``Nonlinear Dispersive
Equations.'' We met for the first time at this program and began
work on this project there.  Also, we are grateful to MSRI for
providing accommodations for S.R. during a May 2006 visit to U.C.\
Berkeley.  Finally, we are indebted to Jim Colliander for mentorship
and encouragement. J.H. is partially supported by an NSF
postdoctoral fellowship.  S.R. is partially supported by NSF grant
DMS-0531337.

\section{Dichotomy for the energy subcritical NLS}
\label{S:finite-time}

In this section we study a more general version of the focusing
nonlinear Schr\"odinger equation NLS$_p(\cR^N)$ which is mass
supercritical and energy subcritical, i.e.
\begin{equation}
\left\{ \begin{array}{l}
i\partial_t u +\Delta u + |u|^{p-1}u=0,
\quad (x,t) \in \cR^N
\times \cR,\\
u(x,0) = u_0(x),
\end{array}
\right.
\end{equation}
with the choice of nonlinearity $p$ and the dimension $N$ such that
$$
0<s_c<1, \quad \text{where} \quad s_c = \frac{N}2-\frac2{p-1}.
$$
In other words, we consider $\dot{H}^{s_c}$-critical NLS equations
with $0<s_c<1$. In this case the initial value problem with $u_0 \in
{H}^1(\cR^N)$ is locally well-posed, see \cite{GV79}. Denote by $I =
(-T_*, T^*)$ the maximal interval of existence of the solution $u$
(e.g., see \cite{Caz-book}). This implies that either $T^* =
+\infty$ or $T^* < +\infty$ and $\Vert \nabla u (t) \Vert_{L^2} \to
\infty$ as $t \to T^*$ (similar properties for $T_*$).

The solutions to this problem satisfy mass and energy conservation
laws
$$
M[u](t) = \int |u(x,t)|^2 \,dx = M[u_0],
$$
$$
E[u](t) = \frac12\int |\nabla u(x,t)|^2 - \frac1{p+1}\int
|u(x,t)|^{p+1} \,dx = E[u_0],
$$
and the Sobolev $\dot{H}^{s_c}$ norm and Lebesgue $L^{p_c}$ norm,
$p_c = \frac{N}2 (p-1)$, are invariant under the scaling $u \mapsto
u_{\lambda}(x,t) = \lambda^{2/(p-1)} u(\lambda x, \lambda^2 t)$.
(Note that $u_\lambda$ is a solution of NLS$_p(\cR^N)$, if $u$ is.)

We investigate other scaling invariant quantities besides the above
norms. Since
$$
\Vert u_\lambda \Vert_{L^2(\cR^N)} = \lambda^{-s_c}
\Vert u \Vert_{L^2(\cR^N)} \quad \text{and} \quad \Vert \nabla
u_\lambda \Vert_{L^2(\cR^N)} = \lambda^{-s_c+1}
\Vert \nabla u \Vert_{L^2(\cR^N)},
$$
the quantity (or any power of it)
$$
\Vert \nabla u_0 \Vert^{s_c}_{L^2(\cR^N)} \cdot \Vert u_0
\Vert_{L^2(\cR^N)}^{1-s_c} %
$$
is scaling invariant. Another scaling invariant quantity
is
$$ \Lambda_0:= E[u_0]^{s_c}\, M[u_0]^{1-s_c}.
$$

Next, recall the Gagliardo-Nirenberg inequality from \cite{W83}
which is valid for values $p$ and $N$ such that $0 \leq s_c <
1$\footnote{It is also valid for $s_c=1$ becoming nothing else but
Sobolev embedding, see Remark \ref{W<->KM}.
}:
\begin{equation}
 \label{GNgeneral}
\Vert u \Vert^{p+1}_{L^{p+1}(\cR^N)} \leq c_{GN} \, \Vert \nabla u
\Vert_{L^2(\cR^N)}^{\frac{N(p-1)}{2}} \, \Vert u
\Vert_{L^2(\cR^N)}^{2-\frac{(N-2)(p-1)}{2}},
\end{equation}
where $\ds c_{GN} = c_{GN}(p,N) = \frac{p+1}{2 \,\Vert Q
\Vert_{L^2(\cR^N)}^{p-1}}$ and $Q$ is the ground state solution
(positive solution of minimal $L^2$ norm) of the following
equation\footnote{We use the notation from Weinstein \cite{W83}; one
can rescale $Q$ so it solves $\lap Q -Q +Q^p = 0$.}
\begin{equation}
 \label{E:Qground}
\frac{N(p-1)}4\, \lap Q - \left(1-\frac{(N-2)(p-1)}4\right)\, Q +
\,Q^{p} = 0.
\end{equation}
(See \cite{W83} and references therein for the discussion on the
existence of positive solutions of class $H^1(\cR^N)$ to this
equation.)

Define
$$
u_Q(x,t)= e^{i\lambda t} Q(\alpha x)
$$
with $\alpha = \frac{\sqrt{N(p-1)}}{2}~ \text{and} ~ \lambda = 1-
\frac{(N-2)(p-1)}4$. Then $u_Q$ is a soliton solution of $i
\partial_t u +\lap u + |u|^{p-1} u = 0$.
Since the Gagliardo-Nirenberg inequality is optimized by $Q$, we get
\begin{equation}
 \label{E:Qground1}
\Vert Q \Vert_{L^{p+1}(\cR^N)}^{p+1}
= \frac{p+1}2 \, \Vert
\nabla Q \Vert_{L^2(\cR^N)}^{\frac{N(p-1)}2} \,\Vert Q
\Vert_{L^2(\cR^N)}^{2-\frac{N(p-1)}2}.
\end{equation}
Multiplying (\ref{E:Qground}) by $Q$ and integrating, we obtain
\begin{equation}
 \label{E:Qground2}
-\lambda \Vert Q \Vert^2_{L^2(\cR^N)} - \alpha^2 \Vert \nabla Q
\Vert^2_{L^2(\cR^N)} + \Vert Q \Vert^{p+1}_{L^{p+1}(\cR^N)} = 0.
\end{equation}
Combining (\ref{E:Qground1}) and (\ref{E:Qground2}), we obtain
$$
\frac{p+1}2 \, \Vert Q \Vert^2_{L^2(\cR^N)} \, \Vert \nabla Q
\Vert_{L^2(\cR^N)}^{\frac{N(p-1)}2} = \lambda \, \Vert Q
\Vert^{2+\frac{N(p-1)}2}_{L^2(\cR^N)} + \alpha^2 \, \Vert \nabla Q
\Vert_{L^2(\cR^N)}^2 \, \Vert Q \Vert^{\frac{N(p-1)}2}_{L^2(\cR^N)}.
$$
Note that $\Vert Q \Vert_{L^2(R^N)} = 0$ is the trivial solution of
the above equation and we exclude it from further consideration.
Denote by $\ds z = \frac{\Vert \nabla Q \Vert_{L^2(R^N)}}{\Vert Q
\Vert_{L^2(R^N)}}$.  Now the above equation becomes
$$
\frac{p+1}2 \, z^{\frac{N(p-1)}2} - \frac{N(p-1)}4 z^2 +
\frac{(N-2)(p-1)}4 - 1 = 0.
$$

The equation has only one real root $z=1$ which gives
\begin{itemize}
\item[(Q.1)]
$\Vert \nabla Q \Vert_{L^2(\cR^N)} = \Vert Q \Vert_{L^2(\cR^N)}$.
\end{itemize}

Substituting (Q.1) into (\ref{E:Qground1}), we also obtain
\begin{itemize}
\item[(Q.2)]
$\Vert Q \Vert^{p+1}_{L^{p+1}(\cR^N)} = \frac{p+1}2 \, \Vert Q
\Vert^2_{L^2(\cR^N)}$.
\end{itemize}

We are now ready to state the main result of this section.
\begin{theorem}
 \label{BlowUp}
Consider NLS$_p(\cR^N)$ with $u_0 \in {H}^1(\cR^N)$ and $0 < s_c <
1$. Let $u_Q (x,t)$ be as above and denote $\ds \sigma_{p,N} =
\left( \frac{4}{N(p-1)}\right)^{\frac1{p-1}} \, \Vert Q
\Vert_{L^2(\cR^N)}$. Suppose that
\begin{equation}
 \label{E:initialLambda}
\Lambda_0 := E[u_0]^{s_c}\, M[u_0]^{1-s_c} < E[u_Q]^{s_c} \,
M[u_Q]^{(1-s_c)} \equiv \left(\frac{s_c}{N}\right)^{s_c}
(\sigma_{p,N})^2, \quad E[u_0] \geq 0.
\end{equation}
If
\begin{equation}
 \label{E:less}
\Vert \nabla u_0 \Vert^{s_c}_{L^2} \cdot \Vert u_0
\Vert_{L^2}^{1-s_c} < \Vert \nabla u_Q \Vert^{s_c}_{L^2} \cdot \Vert
u_Q \Vert_{L^2}^{1-s_c} \equiv \sigma_{p,N},
\end{equation}
then $I = (-\infty, +\infty)$, i.e. the solution exists globally in
time, and for all time $t \in \cR$
\begin{equation}
 \label{E:lessallt}
\Vert \nabla u(t) \Vert^{s_c}_{L^2} \cdot \Vert u_0
\Vert_{L^2}^{1-s_c} < \Vert \nabla u_Q \Vert^{s_c}_{L^2} \cdot \Vert
u_Q \Vert_{L^2}^{1-s_c} \equiv \sigma_{p,N}.
\end{equation}
If
\begin{equation}
 \label{E:greater}
\Vert \nabla u_0 \Vert^{s_c}_{L^2(\cR^N)} \cdot \Vert u_0
\Vert_{L^2(\cR^N)}^{1-s_c} > \Vert \nabla u_Q \Vert^{s_c}_{L^2}
\cdot \Vert u_Q \Vert_{L^2}^{1-s_c} \equiv \sigma_{p,N},
\end{equation}
then for $t \in I$
\begin{equation}
 \label{E:greaterallt}
\Vert \nabla u(t) \Vert^{s_c}_{L^2(\cR^N)} \cdot \Vert u_0
\Vert^{1-s_c}_{L^2(\cR^N)} > \Vert \nabla u_Q \Vert^{s_c}_{L^2}
\cdot \Vert u_Q \Vert_{L^2}^{1-s_c} \equiv \sigma_{p,N}.
\end{equation}
Furthermore, if $|x|u_0 \in L^2(\cR^N)$, then $I$ is finite, and
thus, the solution blows up in finite time. The finite-time blowup
conclusion and \eqref{E:greaterallt} also hold if, in place of
\eqref{E:initialLambda} and \eqref{E:greater}, we assume $E[u_0]<0$.

\end{theorem}

\begin{remark}
It is easy to check the equivalence on the right-hand side of
(\ref{E:initialLambda}) - (\ref{E:greaterallt}):
\begin{align*}
E[u_Q]^{s_c} M[u_Q]^{(1-s_c)} &= \left( \frac{\alpha^2}2 -
\frac12\right)^{s_c} \alpha^{-N}\, \Vert Q \Vert_{L^2(\cR^N)}^2\\
&= \left(\frac{N(p-1)-4}8\right)^{s_c} \, \left(
\frac4{N(p-1)}\right)^{\frac{N}2} \, \Vert Q \Vert_{L^2(\cR^N)}^2.
\end{align*}
On the other hand,
$$
\left( \frac{s_c}{N}\right)^{s_c} (\sigma_{p,N})^2 = \left(
\frac{N(p-1)-4}{2N(p-1)} \right)^{s_c} \, \left(
\frac4{N(p-1)}\right)^{\frac2{(p-1)}} \, \Vert Q
\Vert_{L^2(\cR^N)}^2,
$$
which equals the previous expression when recalling that $s_c
=\frac{N}2 - \frac2{p-1}$.

Furthermore, since $\sigma_{p,N} = \left(\frac4{N(p-1)}
\right)^{\frac1{(p-1)}} \Vert Q \Vert_{L^2(\cR^N)}$, we also obtain
that
$$
\Vert \nabla u_Q \Vert^{s_c}_{L^2(\cR^N)} \cdot \Vert u_Q
\Vert^{(1-s_c)}_{L^2(\cR^N)} = \alpha^{s_c-\frac{N}2} \Vert Q
\Vert_{L^2(\cR^N)} = \alpha^{-\frac2{p-1}} \, \Vert Q
\Vert_{L^2(\cR^N)} \equiv \sigma_{p,N}.
$$
\end{remark}
\smallskip

\begin{remark}
Observe that the second part of Theorem \ref{BlowUp} shows that
there are solutions of NLS$_p(\cR^N)$ with positive energy which
blow up in finite time, thus, we extend the standard virial argument
(e.g. see \cite{G77}) on the existence of blow up solutions with
negative energy and finite variance. Moreover, using the localized
version of the virial identity, this result can be extended to the
functions with infinite variance, see Corollary \ref{NoVariance}
below.
\end{remark}

\begin{remark} \label{W<->KM}
This theorem provides a link between the mass critical NLS and
energy critical NLS equations. Consider $s_c=1$, then the theorem
holds true by the work of Kenig-Merle \cite[Section
3]{KM}: in this case %
$\Lambda_0 = E[u_0]$, the Gagliardo-Nirenberg inequality
(\ref{GNgeneral}) becomes the Sobolev inequality with $c_{N} =
(c_{GN})^{1/(p+1)}$, the condition (\ref{E:initialLambda}) becomes
$E[u_0] < \frac1{N} \,(c_{GN})^{-\frac2{p-1}} = \frac1{N}\,
(c_{N})^{-N} = E[W]$, where $W$ is the radial positive decreasing
(class $\dot{H}^1(\cR^N)$) solution of $\lap W + |W|^{p-1} W = 0$,
and the conditions (\ref{E:less}) - (\ref{E:greaterallt}) involve
only the size of $\Vert \nabla u_0 \Vert_{L^2}$ in correlation with
$\sigma_{p,N} = (c_{GN})^{-\frac1{p-1}} = c_N^{-N/2} = \Vert \nabla
W \Vert_{L^2(\cR^N)}$.

In the case $s_c = 0$, the only relevant scaling invariant quantity
is the mass: $\Lambda_0 = M[u_0]$, the condition
(\ref{E:initialLambda}) becomes $M[u_0] < \left(\frac{p+1}{2}
\frac1{c_{GN}}\right)^{\frac2{p-1}} = \Vert Q \Vert^2_{L^2(\cR^N)}$,
and the conditions (\ref{E:less}) - (\ref{E:greaterallt}) involve
also only the mass in relation with $\Vert Q \Vert^2_{L^2(\cR^N)}$,
in fact (\ref{E:less}) (and (\ref{E:lessallt})) coincides with
(\ref{E:initialLambda}) and the conclusion on the global existence
holds; the condition (\ref{E:greater}) becomes $\Vert u_0
\Vert_{L^2} > \Vert Q \Vert_{L^2}$, and thus, the complement of
(\ref{E:initialLambda}) holds, which is the last statement
$E[u_0]<0$, hence, implying the blow up.
Thus, the statement of the theorem in the limiting case $s_c=0$
connects with Weinstein's results \cite{W83}.
\end{remark}

\begin{proof}[Proof of Theorem \ref{BlowUp}]
Using the definition of energy and (\ref{GNgeneral}), we have
\begin{align*}
E[u] &= \frac12 \,\Vert \nabla u \Vert_{L^{2}}^2 - \frac1{p+1}
\,\Vert u
\Vert^{p+1}_{L^{p+1}} \\
& \geq \frac12 \,\Vert \nabla u \Vert_{L^{2}}^2 -
\frac{c_{GN}}{p+1}\, \Vert \nabla u \Vert_{L^{2}}^{\frac{N(p-1)}2}
\, \Vert u_0 \Vert^{2-\frac{(N-2)(p-1)}2}_{L^2}.
\end{align*}
Define $\ds f(x) = \frac12 \, x^2 - \frac{c_{GN}}{p+1}\, \Vert u_0
\Vert_{L^{2}}^{2-\frac{(N-2)(p-1)}2} \, x^{\frac{N(p-1)}2}$, observe
that $\deg(f) \geq 2$, since $N(p-1) \geq 4$. Then
\begin{align*}
f'(x) &= x - \frac{N(p-1)}{2(p+1)} c_{GN} \, \Vert u_0
\Vert_{L^{2}}^{2-\frac{(N-2)(p-1)}2} \, x^{\frac{N(p-1)}2 -1}\\
& = x\left (1- \frac{N(p-1)}{2(p+1)} c_{GN} \, \Vert u_0
\Vert_{L^{2}}^{(p-1)(1-s_c)} \, x^{(p-1)s_c} \right),
\end{align*}
and thus, $f'(x) = 0$ when $x_0=0$ and $\ds x_1 =
(\sigma_{p,N})^{\frac1{s_c}} \, \Vert u_0
\Vert_{L^2(\cR^N)}^{-\frac{1-s_c}{s_c}}$.
Note that $f(0)=0$ and $\ds f(x_1) = \frac{s_c}{N} \, x_1^2$. Thus,
the graph of $f$ has two extrema: a local minimum at $x_0$ and a
local maximum at $x_1$. Hence, the condition (\ref{E:initialLambda})
implies that $E[u_0] < f(x_1)$.
Combining this with energy conservation, we have
\begin{equation}
\label{E:coercion} f(\|\nabla u(t)\|_{L^2}) \leq
E[u(t)]=E[u_0]<f(x_1).
\end{equation}

If initially $\Vert \nabla u_0 \Vert_{L^2} < x_1$, i.e. the
condition (\ref{E:less}) holds, then by \eqref{E:coercion} and the
continuity of $\|\nabla u(t)\|_{L^2}$ in $t$, we have $\|\nabla u(t)\|_{L^2} <
x_1$ for all time $t \in I$ which gives (\ref{E:lessallt}). In
particular, the $\dot{H}^1$ norm of the solution $u$ is bounded,
which proves global existence in this case.

If initially $\Vert \nabla u_0 \Vert_{L^2} > x_1$, i.e. the
condition (\ref{E:greater}) holds, then by \eqref{E:coercion} and
the continuity of $\|\nabla u(t)\|_{L^2}$ in $t$, we have $\|\nabla u(t)\|_{L^2} >
x_1$ for all time $t \in I$ which gives (\ref{E:greaterallt}).
Now if $u$ has the finite variance, we recall the virial identity
$$
\partial_t^2 \int |x|^2 \, |u(x,t)|^2 \, dx = 4 N(p-1) E[u_0] -
2(N(p-1)-4)\Vert \nabla u (t) \Vert^2_{L^2}.
$$
Multiplying both sides by $M[u_0]^\theta$ with $\theta =
\frac{1-s_c}{s_c}$ and applying inequalities (\ref{E:initialLambda})
and (\ref{E:greaterallt}), we obtain
\begin{align*}
M[u_0]^{\theta} \, \partial_t^2 \int |x|^2 \, |u(x,t)|^2 \, dx & = 4
N(p-1) \Lambda_0^{\frac1{s_c}} - 2(N(p-1)-4)\Vert \nabla u (t)
\Vert^2_{L^2} \,
\Vert u_0 \Vert^{2 \,\theta}_{L^2}\\
&< 4(p-1) s_c (\sigma_{p,N})^{\frac2{s_c}} - 2(N(p-1) - 4)
(\sigma_{p,N})^{\frac2{s_c}} = 0,
\end{align*}
and thus, $I$ must be finite, which implies that in this case blow
up occurs in finite time.\footnote{To be more accurate, in order to
obtain the finite-time blow-up one in fact needs to deduce from
\eqref{E:coercion} that $\|\nabla u(t)\|_{L^2} \geq x_1+\delta_1$
and consequently that $M[u_0]^{\theta} \, \partial_t^2 \int |x|^2 \,
|u(x,t)|^2 \, dx \leq -\delta_2 <0$, where $\delta_1>0$ can be
determined in terms of  $f(x_1)-E[u_0]>0$ and $\delta_2>0$ in terms
of $\delta_1$.}
\end{proof}

\begin{corollary}
 \label{NoVariance}
Suppose that all conditions of Theorem \ref{BlowUp} hold except for
finite variance, and now assume that the solution $u$ is radial.
Consider $N \geq 2$ and $1+\frac4{N} < p < \min\{1+ \frac4{N-2}, 5
\}$\footnote{This is a technical restriction.}. Also suppose that
there exists $\delta>0$ such that
\begin{equation}
 \label{E:initialLambda-delta}
\Lambda_0 \leq (1-\delta)^{s_c} \left(\frac{s_c}{N} \right)^{s_c}
(\sigma_{p,N})^2.
\end{equation}
If (\ref{E:greater}) holds, then there exists $\tilde{\delta} =
\tilde{\delta}(\delta) > 0$ such that
\begin{equation}
 \label{E:greaterallt-delta}
\Vert \nabla u(t) \Vert^{s_c}_{L^2(\cR^N)} \cdot \Vert u_0
\Vert^{1-s_c}_{L^2(\cR^N)} \geq (1+ \tilde\delta)^{s_c}
\,\sigma_{p,N}.
\end{equation}
Furthermore, the maximal interval of existence $I$ is finite.
\end{corollary}

\begin{proof}
The inequality (\ref{E:greaterallt-delta}) follows from the proof of
the theorem by applying a refined version
(\ref{E:initialLambda-delta}), so we concentrate on the second
implication. We use a localized version of the virial identity (e.g.
\cite{KM}), let $\varphi \in C_0^{\infty}(\cR^N)$, then
\begin{align*}
\partial_t^2 \int \varphi(x) \, |u(x,t)|^2 \, dx &=
4  \sum_{j,k} \R \int \partial_{x_j}\partial_{x_k} \varphi \,
\partial_{x_j} u \, \partial_{x_k} \bar{u} - \int \lap^2 \varphi
\, |u|^2\\
&\qquad - 4\left(\frac12 - \frac1{p+1}\right) \int \lap \varphi \,
|u|^{p+1}.
\end{align*}
Choose $\varphi(x)=\varphi(|x|)$ to be a radially symmetric function
that is constant for large $r$ and such that $\partial_r^2
\varphi(r) \leq 2$ for all $r\geq 0$ and $\varphi(r) = r^2$ for
$0\leq r \leq 1$.  Taking $\varphi_m(r) = m^2\varphi(r/m)$, and
following the proof of the main theorem in Ogawa-Tsutsumi
\cite{OT91}\footnote{This proof uses the radial Gagliardo-Nirenberg
estimate, and hence, we have the radial restriction in Cor
\ref{NoVariance}.}, we obtain that for any large $m > 0$ and $\gamma
= (N-1)(p-1)/2$, we have
\begin{align*}
\partial_t^2 \int \varphi_m(|x|) \, |u(x,t)|^2 \, dx &\leq
4N(p-1) E[u_0] - (2N(p-1) - 8) \int |\nabla
u|^2\\
&\qquad + \frac{c_1}{m^{\gamma}} \, \Vert u \Vert_{L^2}^{(p+3)/2} \, \Vert
\nabla u \Vert^{(p-1)/2}_{L^2} + \frac{c_2}{m^2} \, \int_{m<|x|}
|u|^2.
\end{align*}
Choose
$$
0 < \ep < (2N(p-1)-8)\left( 1- \frac{1-\delta}{(1+\tilde\delta)^2}
\right).
$$
Use Young's inequality in the third term on the right
hand side to separate the $L^2$-norm and gradient term and then
absorb the gradient term into the second term with the chosen $\ep$.
Then multiplying the above expression by $M[u_0]^{\theta}$ and using
(\ref{E:initialLambda-delta}) and (\ref{E:greaterallt-delta}), we
get
\begin{align*}
\indentalign M[u_0]^{\theta} \, \partial_t^2 \int \varphi_m(|x|) \, |u(x,t)|^2 \,
dx \\
&\leq 4N(p-1) \, E[u_0] M[u_0]^{\theta} - (2N(p-1) - 8 -\ep ) \Vert \nabla u
\Vert^2_{L^2} \Vert u \Vert^{2\theta}_{L^2}\\
&\qquad + \frac{c(\ep,N)}{m^{4\gamma/(5-p)}} \, \Vert u
\Vert_{L^2}^{2(p+3)/(5-p)+2\theta} + \frac{c_2}{m^2} \, \Vert u
\Vert_{L^2}^{2+2\theta}\\
&\leq 4N(p-1) \, (1-\delta) \, \frac{s_c}{N} \,(\sigma_{p,N})^{\frac2{s_c}} - (2N(p-1) - 8 -\ep ) (1+ \tilde\delta)^2\,
(\sigma_{p,N})^{\frac2{s_c}}\\
&\qquad + \frac{c(\ep,N)}{m^{4\gamma/(5-p)}} \, \Vert u
\Vert_{L^2}^{2(p+3)/(5-p)+2\theta} + \frac{c_2}{m^2} \, \Vert u
\Vert_{L^2}^{2+2\theta}\\
& < - c(\epsilon, N,p)
\end{align*}
by choosing $m=m(\epsilon, \delta, \tilde \delta, N,p, M[u_0])$
large enough, where $c(\epsilon, N,p)>0$. This implies that the
maximal interval of existence $I$ is finite.
\end{proof}

\section{Critical norm concentration phenomenon}
\label{S:conc}

This section is devoted to a proof of Theorem \ref{T:L3conc}. We
will prove Proposition \ref{P:L3conc} below, and then indicate how
Theorem \ref{T:L3conc} is an easy corollary of this proposition.

To state the proposition, we need some notation for spatial and
frequency localizations.  Let $\phi(x)\in C_c^\infty(\mathbb{R}^3)$
be a radial function so that $\phi(x)=1$ for $|x|\leq 1$ and
$\phi(x)=0$ for $|x|\geq 2$, and then define the inner and outer
spatial localizations of $u(x,t)$ at radius $R(t)>0$ as $u_1(x,t) =
\phi(x/R(t))u(x,t)$, $u_2(x,t)=(1-\phi(x/R(t)))u(x,t)$. Let
$\chi(x)\in C_c^\infty(\mathbb{R}^3)$ be a radial function so that
$\chi(x)=1$ for $|x|\leq 1/8\pi$ and $\chi(x)=0$ for $|x|\geq
1/2\pi$, and furthermore, $\hat\chi(0)=1$ and define the inner and
outer frequency localizations at radius $\rho(t)$ of $u_1$ as $\hat
u_{1L}(\xi,t) = \hat\chi(\xi/\rho(t))\hat u_1(\xi,t)$ and $\hat
u_{1H}(\xi,t) = (1-\hat\chi(\xi/\rho(t)))\hat
u_1(\xi,t)$.\footnote{The $1/8\pi$ and $1/2\pi$ radii are chosen to
be consistent with the assumption $\hat \chi(0)=1$, since
$\hat\chi(0)=\int_{\mathbb{R}^3} \chi(x) dx$.  In actuality, this is
for convenience only--the argument is easily adapted to the case
where $\hat\chi(0)$ is any number $\neq 0$.}  Note that the
frequency localization of $u_1=u_{1L}+u_{1H}$ is inexact, although
crucially we have
\begin{equation}
\label{E:smallfreq}
|1-\hat\chi(\xi)|\leq c\min(|\xi|,1) \,.
\end{equation}

\begin{proposition}
\label{P:L3conc} Let $u$ be an $H^1$ radial solution to
\eqref{E:NLS} that blows-up in finite $T>0$. Let $R(t)=
c_1\|u_0\|_{L^2}^{3/2}\|\nabla u(t)\|_{L^2_x}^{-1/2}$ and $\rho(t) =
c_2\|\nabla u(t)\|_{L_x^2}^2$ (for absolute constants $c_1$ and
$c_2$), and decompose $u=u_{1L}+u_{1H}+u_2$ as described in the
paragraph above.
\begin{enumerate}
\item There exists an absolute constant $c>0$ such that
\begin{equation}
\label{E:u1Lbd}
\|u_{1L}(t)\|_{L_x^3} \geq c \quad\text{as}\quad t\to T.
\end{equation}
\item
Suppose that there exists a constant $c^*$ such that
$\|u_1(t)\|_{L^3} \leq c^*$. Then
\begin{equation}
\label{E:narrow} \|u_1(t)\|_{L_x^3(|x-x_0(t)|\leq \rho(t)^{-1})}
\geq \frac{c}{(c^*)^3} \quad\text{as}\quad t\to T
\end{equation}
for some absolute constant $c>0$, where $x_0(t)$ is a position
function such that $|x_0(t)|/{\rho(t)^{-1}} \leq c\cdot (c^*)^6$.
\end{enumerate}
\end{proposition}

We mention two Gagliardo-Nirenberg estimates for functions on
$\mathbb{R}^3$ that will be applied in the proof.  The first is an
``exterior'' estimate, applicable to radially symmetric functions
only, originally due to W. Strauss:
\begin{equation}
 \label{E:GNout}
\| v \|_{L_{\{|x|>R\}}^4}^4 \leq \frac{c}{R^2}
\|v\|_{L_{\{|x|>R\}}^2}^3 ||\nabla v \|_{L_{\{|x|>R\}}^2},
\end{equation}
where $c$ is an absolute constant (in particular, independent of
$R>0$). The second is generally applicable:  For any function $v$,
\begin{equation}
 \label{E:GNin}
\|v\|^4_{L^4(\cR^3)} \leq c\|v\|_{L^3(\cR^3)}^2 \|\nabla
v\|_{L^2(\cR^3)}^2 \,.
\end{equation}

\begin{proof}[Proof of Prop.\ \ref{P:L3conc}]
Since, by \eqref{E:lowerbd}, $\|\nabla u(t)\|_{L_x^2} \to +\infty$
as $t\to T$, by energy conservation, we have
$\|u(t)\|_{L_x^4}^4/\|\nabla u(t)\|_{L_x^2}^2 \to 2$.  Thus, for $t$
sufficiently close to $T$,
\begin{equation}
\label{E:en} \|\nabla u\|_{L_x^2}^2 \leq \|u\|_{L_x^4}^4 \leq
\|u_{1L}\|_{L_x^4}^4 + \|u_{1H}\|_{L_x^4}^4+\|u_2\|_{L_x^4}^4 \,.
\end{equation}
By \eqref{E:GNout}, the choice of $R(t)$, and mass conservation
\begin{equation}
\label{E:u2bd} \|u_2\|_{L_x^4}^4 \leq
\frac{c}{R^2}\|u_0\|_{L_x^2}^3\|\nabla u\|_{L_x^2} \leq
\frac14\|\nabla u\|_{L_x^2}^2,
\end{equation}
where the constant $c_1$ in the definition of $R(t)$ has been chosen
to obtain the factor $\frac14$ here.  By Sobolev embedding,
\eqref{E:smallfreq}, and the choice of $\rho(t)$,
\begin{align}
\|u_{1H}\|_{L_x^4}^4 &\leq c \, \|u_{1H}\|_{\dot H_x^{3/4}}^4
= c \, \| |\xi|^{3/4}(1-\hat\chi(\xi/\rho))\hat u_1(\xi)\|_{L_\xi^2}^4 \notag\\
&\leq c \, \rho^{-1}\| |\xi|\hat u_1(\xi)\|_{L_\xi^2}^4 \leq c \,
\rho^{-1}\|\nabla u_1\|_{L_x^2}^4
\leq c \, \rho^{-1}\|\nabla u\|_{L_x^2}^4 \notag \\
&\leq \frac14\|\nabla u\|_{L_x^2}^2,
 \label{E:u1Hbd}
\end{align}
where the constant $c_2$ in the definition of $\rho(t)$ has been
chosen to obtain the factor $\frac14$ here.  Combining \eqref{E:en},
\eqref{E:u2bd}, and \eqref{E:u1Hbd}, we obtain
\begin{equation}
\label{E:u1L_lowbd} \|\nabla u\|_{L_x^2}^2 \leq
c\|u_{1L}\|_{L_x^4}^4.
\end{equation}
By \eqref{E:u1L_lowbd} and \eqref{E:GNin}, we obtain \eqref{E:u1Lbd}, completing the proof of part (1) of the proposition.  To prove part (2), we assume $\|u_1\|_{L^3} \leq c^*$.  By \eqref{E:u1L_lowbd},
\begin{align*}
\|\nabla u\|_{L^2}^2 &\leq c\|u_{1L}\|_{L_x^3}^3\|u_{1L}\|_{L_x^\infty}
\leq c\cdot (c^*)^3\|u_{1L}\|_{L_x^\infty}\\
&\leq c\cdot (c^*)^3 \sup_{x\in \mathbb{R}^3} \left| \int
\rho^3\chi(\rho(x-y))u_1(y)\, dy \right|.
\end{align*}
There exists $x_0=x_0(t)\in \mathbb{R}^3$ for which at least
$\frac12$ of this supremum is attained.  Thus,
\begin{align*}
\|\nabla u\|_{L^2}^2 &\leq c\cdot (c^*)^3  \left| \int
\rho^3\chi(\rho(x_0-y))u_1(y)\,
dy \right|\\
&\leq c\cdot (c^*)^3 \rho^3 \int_{|x_0(t)-y|\leq \rho^{-1}}|u_1(y)|\, dy\\
&\leq c\cdot (c^*)^3 \rho \left( \int_{|x_0(t)-y|\leq \rho^{-1}}
|u_1(y)|^3 \, dy \right)^{1/3}
\end{align*}
with H\"older's inequality used in the last step. By the choice of
$\rho$, we obtain \eqref{E:narrow}. To complete the proof, it
remains to obtain the stated control on $x_0(t)$, which will be a
consequence of the radial assumption and the assumed bound
$\|u_1\|_{L^3}\leq c^*$. Suppose
$$
\frac{|x_0(t_n)|}{\rho(t_n)^{-1}} \gg (c^*)^6
$$
along a sequence of times $t_n\to T$.  Consider the spherical annulus
$$
A = \{ \, x\in \mathbb{R}^3\, : \, |x_0|-\rho^{-1} \leq |x| \leq
|x_0|+\rho^{-1}\, \}
$$
and inside $A$ place $\sim \frac{4\pi |x_0|^2}{\pi (\rho^{-1})^2}$
disjoint balls, each of radius $\rho^{-1}$, centered on the sphere
at radius $|x_0|$.  By the radiality assumption, on each ball $B$,
we have $\|u_1\|_{L_B^3} \geq c/(c^*)^3$, and hence on the annulus
$A$,
$$
\|u_1\|_{L_A^3}^3 \geq
\frac{c}{(c^*)^9}\frac{|x_0|^2}{(\rho^{-1})^2}\gg (c^*)^3,
$$
which contradicts the assumption $\|u_1\|_{L^3} \leq c^*$.
\end{proof}

Now we indicate how to obtain Theorem \ref{T:L3conc} as a consequence.

\begin{proof}[Proof of Theorem \ref{T:L3conc}]
By part (1) of Prop.\ \ref{P:L3conc} and the standard convolution
inequality:
$$
c \leq \|u_{1L}\|_{L_x^3} = \|\rho^3\chi(\rho \,\cdot ) \ast u_1
\|_{L_x^3} \leq \|u_1\|_{L^3}.
$$
Now, if $\|u_1(t)\|_{L^3}$ is not bounded above, then there exists a
sequence of times $t_n\to T$ such that $\|u_1(t_n)\|_{L^3} \to
+\infty$.  Since $\|u(t_n)\|_{L^3(|x|\leq 2R)} \geq
\|u_1(t_n)\|_{L^3}$, we have \eqref{E:wide} in Theorem
\ref{T:L3conc}.  If, on the other hand, $\|u_1(t)\|_{L^3} \leq c^*$,
for some $c^*$, as $t\to T$, we have \eqref{E:narrow} of Prop.
\ref{P:L3conc}.  Since $|x_0(t)|\leq c(c^*)^6\rho(t)^{-1}$, we have
$$
\frac{c}{(c^*)^3} \leq \|u_1(t)\|_{L^3(|x-x_0(t)|\leq \rho(t)^{-1})}
\leq \|u_1(t_n)\|_{L^3(|x|\leq c(c^*)^6\rho(t)^{-1})},
$$
which gives
\eqref{E:tight} in Theorem \ref{T:L3conc}.
\end{proof}

\section{Heuristic analysis of contracting sphere blow-up solutions}
\label{S:heuristic}

In this section, we develop a heuristic analysis of hypothetical
contracting sphere solutions of \eqref{E:NLS}  with blow-up time
$T>0$.  The results of this analysis are summarized in the
introduction \S\ref{S:intro} and we now present the argument itself.
We assume radial symmetry of the solution and let $r=|x|$ denote the
radial coordinate. In radial coordinates, \eqref{E:NLS} becomes
\begin{equation}
\label{E:NLSrad} i\partial_t u + \frac2r\partial_r u + \partial_r^2
u + |u|^2u = 0.
\end{equation}
We define the radial position $r_0(t)$ as
\begin{equation}
\label{E:r0def} r_0^2(t) = \frac{\int_0^{+\infty} r^4|u(r,t)|^2 \,
dr}{\int_0^{+\infty} r^2 |u(r,t)|^2 \, dr} =
\frac{V[u](t)}{M[u](t)},
\end{equation}
where $V[u](t)$ is the first virial quantity
$$
V[u](t) = 4\pi \int_0^{+\infty} r^4|u(r,t)|^2 \,dr,
$$
and $M[u](t)=M[u_0]$ is the conserved mass of $u$. In terms of
$r_0(t)$, we define the focusing factor $\lambda(t)$ via the
relation
\begin{equation}
\label{E:lambdadef} \frac{r_0^2(t)}{\lambda^3(t)} = \|\nabla
u(t)\|_{L^2}^2.
\end{equation}
The motivation for these definitions is that if there were some
fixed profile $\phi(y)$ (where $y\in \mathbb{R}$) and $u(r,t)$ had
approximately the form
$$
|u(r,t)| \approx \frac{1}{\lambda_1(t)} \, \phi\left(
\frac{r-r_1(t)}{\lambda(t)}\right),
$$
then the formulas \eqref{E:r0def} and \eqref{E:lambdadef} would give
$r_0(t)\approx r_1(t)$ and $\lambda(t)\approx\lambda_1(t)$.  The
convenience of these definitions is first that they are
\textit{always} defined (on the maximal existence interval $[0,T)$),
and thus, we can distinguish between contracting sphere and
non-contracting sphere solutions in terms of their behavior, and
second that they immediately relate the parameters $r_0(t)$ and
$\lambda(t)$ to quantities appearing in conservation laws for
\eqref{E:NLS}.

Let $s(t)$ denote the rescaled time
$$
s(t) = \int_0^t \frac{d\sigma}{\lambda(\sigma)^2}
$$
and take $S=s(T)$ to be the rescaled blow-up time (which could, and
in fact, will be shown to be $+\infty$). Since $s(t)$ is
monotonically increasing, we can also work with $t(s)$. Then
$r_0(t)$ and $\lambda(t)$ become functions of $s$ and thus we can
speak of $\lambda(s)\equiv\lambda(t(s))$ and $r_0(s)\equiv
r_0(t(s))$. Denote by $y$ the rescaled radial position
$$
y(r,t) = \frac{r-r_0(t)}{\lambda(t)}.
$$
We now introduce the one-dimensional auxiliary quantity $w(y,s)$,
defined by
\begin{equation}
\label{E:u-w} u(r,t) = \frac{1}{\lambda(t)} \,
w\left(\frac{r-r_0(t)}{\lambda(t)},s(t)\right),
\end{equation}
or equivalently,
$$
w(y,s)=\lambda(s)\, u(\lambda(s)r+r_0(s),t(s)).
$$
It is checked by direct computation that $w$ so defined solves an
NLS-type equation
\begin{equation}
\label{E:w_equation} i\partial_s w + \frac{2\lambda}{r}\partial_y w
+ \partial_y^2 w - i \frac{\lambda_s}{\lambda}\Lambda w - i
\frac{(r_0)_s}{\lambda}\partial_yw + |w|^2 w=0,
\end{equation}
where $\Lambda w = (1+y\partial_y)w$. We expect the term
$-i\frac{(r_0)_s}{\lambda}\partial_yw$ to play a more significant
role in our computation than the analogous term in Rapha\"el's
two-dimensional quintic result \cite{R06}.
Now, with these definitions we make \\

\noindent\textbf{Spherical contraction assumption (SCA)}. Suppose that
$\lambda(t)/r_0(t)\to 0$, and that $w(y,s)$ remains well-localized near $y=0$
for all $s$ as $s\to S$.\\

SCA has, in particular, the effect of driving $y_L\equiv
-r_0(t)/\lambda(t)$, which corresponds to the position $r=0$ in the
original coordinates, to $-\infty$ as $t\to T$ (or $s\to S$). By
SCA, the term $\frac{2\lambda}{r}\partial_yw$ in
\eqref{E:w_equation} should be negligible as $s\to S$. We now make a
second simplifying assumption that we later confirm
is consistent.\\

\noindent\textbf{Focusing rate assumption (FRA)}. Suppose that
$\lambda_s/\lambda \to 0$. Note that $\lambda_s/\lambda = \lambda_t\lambda$,
and if $\lambda(t)$ is given by a power-type expression $\lambda(t)=(T-t)^\alpha$,
then since $\lambda_t\lambda=\alpha(T-t)^{2\alpha-1}$, this assumption is valid
if $\alpha>\frac12$.\\

We have no solid justification for introducing FRA, only that it
makes the analysis more tractable. We shall show that FRA gives rise
to one self-consistent scenario (see the comment at the end of
\S\ref{S:appcons}).\footnote{Locating plausible scenarios is really
the objective here anyhow -- an analysis addressing all
possibilities seems too ambitious at this point.} By FRA, we have
license to neglect the term $-i\frac{\lambda_s}{\lambda}\Lambda w$
in \eqref{E:w_equation}.  This leads to the following simplified
equation in $w(y,s)$
\begin{equation}
 \label{E:simplew_equation}
i\partial_s w + \partial_y^2 w -i \frac{(r_0)_s}{\lambda}\partial_y
w + |w|^2w =0.
\end{equation}

\subsection{Asymptotically conserved quantities for $w$}

Assuming $w$ solves \eqref{E:simplew_equation}, we derive
``asymptotically conserved'' quantities for $w$.  The conservation is
only approximate, since \eqref{E:simplew_equation} is only
approximate.  Moreover, in the calculations, we will routinely
ignore the boundary term $y_L=-r_0(s)/\lambda(s)$ corresponding to
$r=0$ in the integration by parts computations.  (An interpretation of the following computations is that ``asymptotically conserved'' quantities like $M[w](t)$ converge to specific values as $t\to T$, i.e.\ behave like $M[w](t)=M[w](T) + \mathcal{O}(T-t)^\mu$ for some $\mu>0$.)

Define the following quantities (here, $y$ is one-dimensional):
\begin{align*}
&M[w] = \int |w|^2 \, dy,\\
&P[w] = \I \int w\,\overline{\partial_y w} \, dy ,\\
&E[w] = \frac12 \int|\partial_yw|^2 \, dy - \frac14\int |w|^4 \, dy,
\end{align*}
the mass, momentum, and energy of $w$. Now we show that the mass and
momentum of $w$ are asymptotically conserved and the energy satisfies
an \textit{a priori} time dependent equation.

Pair the equation (\ref{E:simplew_equation}) with $\bar w$ and take
$2\times$ the imaginary part:
$$
\partial_s \int |w|^2 \, dy = \frac{2(r_0)_s}{\lambda} \R \int
\partial_y w \, \bar w \, dy=0.
$$
Thus, we have $M[w]$ is approximately constant. By direct substitution of the
equation we have
$$\partial_s P[w] = 0,$$
and thus $P[w]$ is approximately constant.

The next step is to study the energy of $w$. For this we introduce
\begin{equation}
 \label{E:w-tilde}
\tilde w(y,s) = e^{-i\frac{(r_0)_s}{2\lambda}y}w(y,s)
\end{equation}
and from \eqref{E:simplew_equation} we see that $\tilde w$ solves
$$
i\partial_s \tilde w - \left( \frac{(r_0)_s}{2\lambda} \right)_s y
\tilde w + \left( \frac{(r_0)_s}{2\lambda}\right)^2 \tilde w +
\partial_y^2 \tilde w + |\tilde w|^2\tilde w=0.
$$
Pair this equation with $\overline{\partial_s\tilde w}$ and take the
real part to get that $E[\tilde w]$ is approximately constant.
Substituting the definition of $\tilde w$, we get that the quantity
\begin{equation}
 \label{E:energy}
\frac{1}{2}\left( \frac{(r_0)_s}{2\lambda}\right)^2 M[w] +
\frac{1}{2}\frac{(r_0)_s}{\lambda}P[w]+E[w] = E[\tilde w] =
\text{approx.\ const.}
\end{equation}
Note that this only says that $E[w]$ is constant if we were to know
that $(r_0)_s/\lambda$ is constant, a point we discuss next.

\subsection{Consequences of the asymptotic conservation laws}
\label{S:appcons}

Now we work out a consequence of the mass conservation of $u$ and
the mass conservation of $w$.
\begin{align*}
M[u] &= 4\pi\int_0^\infty |u|^2 \, r^2dr\\
&= \frac{4\pi}{\lambda(t)^2} \int \left|w \left(
\frac{r-r_0(t)}{\lambda(t)}\right) \right|^2 r^2 dr.
\end{align*}
Now by SCA, we have
\begin{equation}
 \label{E:mass_conserv}
M[u] \approx \frac{4\pi r_0^2(t)}{\lambda(t)} \int |w(y)|^2 \, dy =
\frac{4\pi r_0^2(t)}{\lambda(t)}M[w].
\end{equation}
An immediate consequence of this is that
\begin{equation}
 \label{E:approximate-rate}
r_0(t) \approx \lambda(t)^{1/2}.
\end{equation}

Next we study the consequence of energy conservation of both the
initial solution $u$ and the rescaled version $w$.
\begin{align*}
E[u] &= \frac12 \int |\nabla u|^2 \, dx
-\frac14 \int |u|^4 \, dx = 4\pi \left(\frac12 \int |\partial_r u|^2 r^2
\, dr - \frac14 \int |u|^4 r^2
dr\right)\\
& = \frac{4\pi}{\lambda^3(t)} \left( \frac12 \int |\partial_y w|^2
r^2 \, dy - \frac14 \int |w|^4 r^2 \,dy \right) \approx 4\pi \,
\frac{r_0^2(t)}{\lambda^3(t)} \, E[w].
\end{align*}
Since $r_0(t)$ contracts at a slower rate than $\lambda(t)$, we
have $\ds \frac{r_0^2(t)}{\lambda^3(t)} \to \infty$ as $t \to T$,
and hence, we must have $E[w](t) \to 0$.

Consider the equation (\ref{E:energy}) again: the condition $E[w](t)
\to 0$ together with $M[w]$, $P[w]$ and $E[\tilde w]$ being constant
forces either $\frac{(r_0)_s}{\lambda}$ to be a constant, denote it
by $\kappa$, or $\frac{(r_0)_s}{\lambda} \to 0$ as $t \to T$. If
$\lambda(t) \sim (T-t)^{\gamma}$, then by (\ref{E:approximate-rate})
we have $r_0(t) \sim (T-t)^{\gamma/2}$, and thus, $\kappa = (r_0)_t
\lambda = c \, (T-t)^{\gamma/2-1+\gamma}$ implies $\gamma = 2/3$ if
$\kappa \neq 0$ and $\gamma > 2/3$ if $\kappa = 0$. The latter case
is ruled out by the virial identity in the next section. Thus,
\begin{equation}
 \label{E:preliminary-rates}
r_0(t) \sim (T-t)^{1/3} \quad \mbox{and} \quad \lambda(t) \sim
(T-t)^{2/3}.
\end{equation}
Note that under these conditions and SCA the second term in
(\ref{E:w_equation}) has the coefficient $\frac{\lambda}{r} \approx
\frac{\lambda}{r_0} \sim (T-t)^{1/3}$, and thus, the decision to
drop it in the analysis close to the blow-up time was, at least,
self-consistent.  Similarly, the fourth term coefficient in
(\ref{E:w_equation}) $\frac{\lambda_s}{\lambda} = \lambda_t \lambda
\sim (T-t)^{1/3}$ and becomes negligible near the blow-up time as
well.

\subsection{Application of the virial identities} \label{S:virial}

By (\ref{E:u-w}) and SCA,
$$
\|\nabla u(t)\|_{L^2}^2 = 4\pi \int |\partial_r u|^2 r^2 \, dr
\approx 4\pi \frac{r^2_0(t)}{\lambda^3(t)} \int |\partial_y w|^2 \,
dy.
$$
Substituting \eqref{E:lambdadef}, we obtain $\Vert \partial_y
w\Vert^2_{L^2} \approx \frac1{4\pi}$. The virial identities are
\begin{equation}
 \label{E:virial1}
\partial_t \int r^4 |u(r,t)|^2 \, dr = 4 \I \int r^3 \, \bar u \,
\partial_r u \, dr
\end{equation}
and
\begin{equation}
 \label{E:virial2}
4\pi \partial_t^2 \int r^4 |u(r,t)|^2 \, dr = 24E[u] - 4\|\nabla
u\|_{L_{xyz}^2}^2.
\end{equation}
The equation \eqref{E:virial2} produces the approximate relation
\begin{equation}
 \label{E:virial1a}
\partial_t^2(r_0^2(t))=\frac{1}{M[u]}\left(24E[u] -
\frac{4r_0^2(t)}{\lambda^3(t)}\right).
\end{equation}
Observe that $\gamma > 2/3$ would
contradict (\ref{E:virial1a}); similarly we cannot have lower order
corrections in $(T-t)^\gamma$ (e.g. $(T-t)^{1/3}
\log^{\gamma_1}(T-t)$). We write
\begin{equation}
 \label{E:r0-lambda}
r_0(t)=\alpha \, (T-t)^{1/3}, \qquad \lambda(t)=\beta (T-t)^{2/3}
\end{equation}
with, as yet, undetermined coefficients $\alpha$, $\beta$. The
relation \eqref{E:virial1a} forces one relation
$$
\beta = \left( \frac{18}{M[u]}\right)^{1/3}.
$$
To pursue this further, we incorporate the quantities $M[w]$,
$P[w]$, $E[w]$, and $E[\tilde w]$ into the analysis.
The first of the virial relations \eqref{E:virial1} gives
$$
(r_0^2(t))' \approx
-\frac{16\pi}{M[u]}\frac{r_0^3(t)}{\lambda^2(t)}P[w],
$$
which produces the relation
\begin{equation}
\label{E:w_momentum} P[w] =
\frac{M[u]}{24\pi}\frac{\beta^2}{\alpha} = \frac{(12 M[u])^{1/3}} {8
\pi \alpha}.
\end{equation}
The mass conservation from (\ref{E:mass_conserv}) gives
\begin{equation}
\label{E:w_mass} M[w] = \frac{M[u]}{4\pi}\frac{\beta}{\alpha^2} =
\frac{ 18^{1/3} \, M[u]^{2/3}}{4\pi \alpha^2}.
\end{equation}
By the energy conservation from the previous section we have
$$
E[w] = \frac{\lambda^3(t)}{r_0^2(t)} \, E[u] =
\frac{\beta^2}{\alpha^2}(T-t)^{4/3} E[u] \to 0.
$$
We now have three quantities: $\alpha$, $M[w]$, $P[w]$, and
two equations \eqref{E:w_momentum} and \eqref{E:w_mass}. We
substitute these values for $P[w]$, $M[w]$ and $E[w]$ into
(\ref{E:energy}) and obtain
$$
E[\tilde w]= - \frac1{16\pi}.
$$
Observe that we are still free to choose $\alpha$ as it does not
affect any of the conservation properties -- this flexibility will
be used in the next section.

\subsection{Asymptotic profile}

Recall from (\ref{E:w-tilde}) and (\ref{E:r0-lambda}) that
$$
\tilde w (y,s) = e^{-i \kappa y/2}\, w(y,s) \quad
\text{with} \quad \kappa = -\alpha \beta/3
$$
satisfies
$$
i\partial_s \tilde w + \left( \frac{\kappa}{2}\right)^2 \tilde w +
\partial_y^2 \tilde w + |\tilde w|^2\tilde w=0.
$$

On the grounds that $\tilde w(y,s)$ is approximately a
global-in-time solution to the one-dimensional cubic NLS, and is
well localized at the origin, the only reasonable asymptotic
configuration is a stationary soliton (see Zakharov-Shabat
\cite{ZS72})\footnote{Also possible are the envelope solitons or
``breathers'' solutions described in \cite{ZS72}, although we choose
not to investigate this possibility here since they are unstable and
the small corrections to the $w$ equation would likely cause them to
collapse to decoupled solitons moving away from each other.}.  Thus,
we assume that as $s\to +\infty$,
\begin{equation}
\label{E:soliton} \tilde w (y,s) = e^{i\theta_0} e^{i \nu s}
P(y+y_0)
\end{equation}
for some fixed phase shift $\theta_0$ and spatial shift $y_0$ (since
$y_0\neq 0$ amounts to a lower-order modification in $r_0(t)$, we
might as well drop it and take $y_0=0$) and where $\nu$ is to be
chosen later.   Then $P$ satisfies
\begin{equation}
 \label{E:P}
-\sigma P + P'' + |P|^2 P = 0, \quad \sigma = \nu -
\frac{\kappa^2}{4}.
\end{equation}
The solution of this equation is
$$
P(y) = e^{i\theta} \sqrt{2 \sigma} \sech (\sqrt{\sigma} y),
$$
here, $\theta$ is arbitrary. By \eqref{E:soliton} and $y_0=0$, we
have $w(y,s) \approx e^{i\theta}e^{i\nu s} e^{i\kappa y/2}P(y)$. The
analysis from the previous section gave
\begin{align}
\label{E:wtilde_energy} &E[\tilde w] = -\frac{1}{16\pi},\\
\label{E:w_energy} &E[w] = 0, \\
\label{E:w_mass2} &M[w]=M[\tilde w] = \frac{18^{1/3}M[u]^{2/3}}{4\pi\alpha^2},\\
\label{E:w_momentum2} &P[w] = \frac{(12M[u])^{1/3}}{8\pi \alpha}.
\end{align}

We choose $\nu$ such that $E[e^{i \kappa y/2} P]=0$. This implies
that $\frac18\kappa^2M[P]+E[P] = 0$ and using $E[P] = -\frac23
\sigma \sqrt \sigma$ and $M[P] = 4 \sqrt \sigma$, we obtain
$$
\sigma = \frac{3}{4}\kappa^2, \quad \text{and hence,} \quad \nu =
\kappa^2.
$$
The equations \eqref{E:wtilde_energy} and $E[P]=-\frac23\sigma^{3/2}$ give
$$
\sigma = \left( \frac{3}{32\pi} \right)^{2/3}.
$$
The equations \eqref{E:w_mass2} and $M[P]=4\sigma^{1/2}$ give
$$
\alpha = \frac{3^{1/6}}{2\pi^{1/3}}M[u]^{1/3}.
$$
Now we conclude with two consistency checks:  The values of $\kappa$
obtained together with the definition $\kappa = -\alpha\beta/3$, and
the formula for $\beta$, are consistent with the value of $\alpha$
obtained here. The value in \eqref{E:w_momentum2} and $P[e^{i\kappa
y/2}P(y)] = -2 \kappa \sigma^{1/2}$ is consistent with the value of
$\alpha$ obtained here.

Pulling all of this information together, we obtain the description
given in the introduction.

\section{Consistency with higher precision computations}
\label{S:refined}

As a consistency check, we show that the result obtained in the
previous section regarding the ``approximate conservation'' of the
mass, momentum, and energy of $w$ stands up to a second-level of
precision.  To do this, we consider (\ref{E:w_equation}) with only
the approximation $r\approx r_0$ in the second term and we leave the
fourth term as is (in the previous section, we completely dropped
the second and fourth terms):
\begin{equation}
 \label{E:w_equation2}
i\partial_s w + \frac{2\lambda}{r_0}\partial_y w + \partial_y^2 w
-i\frac{\lambda_s}{\lambda}\Lambda w - i
\frac{(r_0)_s}{\lambda}\partial_yw + |w|^2 w=0.
\end{equation}

Pairing this equation with $\bar w$ and taking the imaginary part,
we have
$$
\frac12 \, \partial_s \int |w|^2 + \frac{2\lambda}{r_0} \I \int
\partial_y w \, \bar w - \frac{\lambda_s}{\lambda} \R \int \Lambda w
\, \bar w = 0,
$$
which simplifies to
$$
\frac12 \, \partial_s \int |w|^2 - \frac{2\lambda}{r_0} \, P[w] -
\frac{\lambda_s}{2\lambda} \, M[w] = 0.
$$
Applying (\ref{E:r0-lambda}), (\ref{E:w_momentum}) and
(\ref{E:w_mass}) to the second and third terms (note that they
cancel each other out), we obtain the conservation of mass $M[w]$.

To calculate the refined momentum, substitute (\ref{E:w_equation2})
into the definition of $P[w]$ to get
$$
\partial_s P[w] - \frac{4\lambda}{r_0} \int
|\partial_y w|^2 - \frac{2\lambda_s}{\lambda} P[w] = 0,
$$
from which we obtain
$$
\partial_s P[w] = \frac{\beta}{\alpha \pi} (T-t)^{1/3} - \frac{\beta^4}{18 \pi
\alpha} M[u] (T-t)^{1/3} = 0,
$$
and thus, the momentum of $w$ is also preserved.

Note that in this more precise approximation we obtain the
cancelation of the ``error" terms in (\ref{E:w_equation2}) as it was
claimed in the introduction.

The calculation of the refined energy from (\ref{E:w_equation2})
doesn't produce any similar cancelation, however, it confirms
(\ref{E:energy}). We outline it next: first substitute
$$
w(y,s) = e^{i\frac{(r_0)_s}{2\lambda} y}
\tilde{w}(y,s)
$$
to remove $\ds i \frac{(r_0)_s}{\lambda}\partial_y w$ term and
obtain
$$
i \partial_s \tilde w + i \left(\frac{(r_0)_s}{\lambda}\right)_s
\tilde{w} + \frac{2\lambda}{r_0} \,\partial_y\tilde w + i
\frac{(r_0)_s}{r_0} \, \tilde w + \partial^2_y \tilde w  +
\left(\frac{(r_0)_s}{2\lambda}\right)^2 \tilde w - i
\frac{\lambda_s}{\lambda} \, \Lambda \tilde w
$$
$$
+ \frac{\lambda_s}{\lambda} \frac{(r_0)_s}{2\lambda} \, y \tilde w
+|\tilde w|^{\,2} \tilde w = 0.
$$
Next we substitute
$$
\tilde w(y,s) = e^{i \frac{\lambda_s}{4\lambda} y^2} v(y,s)
$$
into the previous equation which results in
\begin{equation}
 \label{E:energy-v}
i \partial_s v + \frac{2\lambda}{r_0}\partial_y v +\partial_y^2 v +
v \left( \left(\frac{(r_0)_s}{2\lambda}\right)^2 + i \left[
\left(\frac{(r_0)_s}{\lambda}\right)_s + \frac{(r_0)_s}{r_0} -
\frac{\lambda_s}{2\lambda} \right]\right)
\end{equation}
$$
+ y v \left(\frac{\lambda_s}{\lambda} \frac{(r_0)_s}{2\lambda} +
\frac{\lambda_s}{r_0}\right)
+ y^2 v \left(\frac{\lambda_s}{2\lambda}\right)^2 + |v|^2 v =0.
$$
We examine coefficients in front of $v$, $y v$ and $y^2 v$. Observe
that $\left(\frac{(r_0)_s}{2\lambda}\right)^2$ is the largest
coefficient by (\ref{E:r0-lambda}) and in fact, is a constant; all
other coefficients have the order of $(T-t)^{1/3}$ or $(T-t)^{2/3}$,
and therefore, we drop corresponding terms from further analysis.
Thus, we obtain
$$
i \partial_s v + \frac{2\lambda}{r_0} \, \partial_y v + \partial_y^2
v + \left(\frac{(r_0)_s}{2\lambda}\right)^2 v  + |v|^2 v \approx 0.
$$
We proceed further and make the substitution
$$
v(y,s) = e^{-\frac{\lambda}{r_0} y} \, \tilde{v}(y,s)
$$
to remove the $\frac{2\lambda}{r_0} \partial_y v$ term in
(\ref{E:energy-v}) and obtain
$$
i \partial_s \tilde v  + \partial_y^2 \tilde v +
\left[\left(\frac{(r_0)_s}{2\lambda}\right)^2 -
\frac{\lambda^2}{r_0^2}\right] \tilde v + |\tilde v|^2 \tilde v \,
e^{-2\frac{\lambda}{r_0} y} \approx 0.
$$
Again finding the order of the coefficients and expanding the
exponent in the nonlinear term, we drop $\frac{\lambda^2}{r_0^2}$ as
well as all positive powers of $-2\frac{\lambda}{r_0}y$ to get
$$
i \partial_s \tilde v  + \partial_y^2 \tilde v +
\left(\frac{(r_0)_s}{2\lambda}\right)^2 \tilde v + |\tilde v|^2
\tilde v \approx 0.
$$
This produces $E[\tilde v] \approx \text{const} + O((T-t)^\gamma)$,
$\gamma>0$. Revealing all substitutions, we obtain $E[\tilde w]
\approx \text{const} + O((T-t)^{\gamma_1})$, $\gamma_1>0$, and
expressing the last approximate identity in terms of $w$ we end up
with (\ref{E:energy}).

\section{General case NLS$_p(\cR^N)$}
\label{S:general}

Consider the mass supercritical focusing NLS$_p(\cR^N)$ equation
with $p> 1+ \frac4{N}$ nonlinearity

\begin{equation}
 \label{E:NLS-p}
i\partial_t u + \Delta u + |u|^{p-1}u = 0,
\end{equation}
for $(x,t) \in \cR^N \times \cR$ with Schwartz class initial-data $u_0\in \mathcal{S}(\mathbb{R}^N)$.
Then the following scaling of the solution is itself a solution:
$$
u_{\lambda}(x,t) = \lambda^{2/(p-1)} u(\lambda x, \lambda^2 t).
$$
The scale invariant Sobolev norm is $\dot H^{s_c}$, where $s_c =
\frac{N}{2}-\frac{2}{p-1}$ (since $p>1+\frac4N$, we have $s_c>0$).
If $s_c> 1$ (prototypical case $N=3$ and $p=7$), then we do not have
a local theory in $H^1$. \footnote{Indeed, it was observed by
Birnir-Kenig-Ponce-Svanstedt-Vega \cite{BKPSV} that one can take a
finite-time radial $H^1$-blow-up solution (whose existence is
guaranteed by the virial identity) and suitably rescale it to obtain
a solution initially arbitrarily small in $H^1$ that blows-up in an
arbitrarily short interval of time.} One does, however, have local
well-posedness in $H^s$ for $s>s_c$ by the Strichartz estimates on a
maximal time interval $[0,T_s)$ with $\lim_{t\to T_s} \|u(t)\|_{H^s}
=+\infty$.  It would appear that if $s_1>s_2>s_c$, then it might be
possible for $T_{s_1}<T_{s_2}$.  However, a \textit{persistence of
regularity} result also follows from the Strichartz estimates and
gives that $T_{s_1}=T_{s_2}$.  Thus, even though an $H^1$ local
theory is absent, there is a clear distinction between global
solutions and finite-time blow-up solutions, and it still makes
sense to speak of ``the blow-up time.''

A heuristic similar to the one presented in \S \ref{S:heuristic} for
a radial blow-up solution of \eqref{E:NLS-p} results in the
following estimation of parameters.
Let
$$
u(r,t) = \frac1{\lambda(t)^{2/(p-1)}} \,
w\left(\frac{r-r_0(t)}{\lambda(t)}, t \right).
$$
Using the conservation of mass as in \S \ref{S:appcons}, we obtain
$$
M[u] \approx \frac{r_0(t)^{(N-1)}}{\lambda(t)^{\frac{5-p}{p-1}}}
\,|\, \mathbb S^{N-1}| \, M[w],
$$
and thus, $r_0(t) \sim \lambda(t)^{\frac{5-p}{(p-1)(N-1)}}$. This
means that for all quintic nonlinearity mass supercritical problems
the radial solution would blow up on a constant radius sphere as in
Rapha\"el's construction \cite{R06}\footnote{Note that this does not
contradict Remark 5.20 in \cite{KM} where it is remarked that for $N\geq 4$
if $\ds \Vert \nabla u(t)\Vert_{L^2_x} \to \infty$ as $t \to T$,
then the blow up necessarily occurs at the origin.
Together with our analysis, it suggests that blow up can occur simultaneously
at the origin and on the sphere of constant radius.}; for $p=7$, the solution would
blow up on a sphere with radius $r_0(t) \nearrow \infty$, i.e. on an
\textit{expanding} sphere.

Using the virial identities as in \S \ref{S:virial}, we obtain
$$
\lambda(t) \sim (T-t)^{\gamma} \quad \text{with} \quad \gamma =
\frac{(p-1)(N-1)}{(p-1)(N-1) + 5-p},
$$
and correspondingly,
$$
r_0(t) \sim (T-t)^{\frac{5-p}{(p-1)(N-1)+ (5-p)}}.
$$

\end{document}